\documentclass[a4paper,12pt]{article}
\usepackage[latin5]{inputenc}
\usepackage[centertags]{amsmath}
\usepackage{amssymb,latexsym}
\usepackage{amsfonts}
\usepackage{layout}
\usepackage{enumerate}
\usepackage{graphics}
\usepackage{graphicx}
\usepackage{fancyhdr}
\usepackage{stmaryrd}
\usepackage{psfrag}
\usepackage[T1]{fontenc}
\usepackage{tikz}
\newtheorem{theorem}{Theorem}[section]

\newtheorem{example}{Example}[section]

\newtheorem{corollary}{Corollary}[section]

\numberwithin{equation}{section}
\title{\bf On Harmonicity and Differential Equations of a Bertrand Curve}
\author{S\"{u}leyman \c{S}enyurt$^1$ and Osman \c{C}ak{\i}r$^{2}$\thanks{corresponding author: osmancakir75@hotmail.com} \\
\small$^{1}$Faculty of Arts and Sciences, Department of Mathematics, \\
\small Ordu University, 52200, Ordu/Turkey\\
\small$^{2}$ Institute of Science, Department of Mathematics, \\
\small Ordu University, 52200, Ordu/Turkey\\
\small $^{1}$senyurtsuleyman52@gmail.com}
\date{}\date{}
\begin{document}
\maketitle

\noindent \textbf{Abstract.}
In the present paper, we give new Frenet formulas for the Bertrand partner curve by taking the advantage of relations between curvatures and a curve itself. Then making use of these formulas we write the differential equations and sufficient conditions of harmonicity of the Bertrand partner curve in terms of the main curve. Finally we exemplify our assertions on the curve helix to see how the formulas we developed work.\\

biharmonic,\, Laplace operator.

\small{\section{Introduction and Preliminaries}}\label{S:intro}
\noindent In differential geometry one of the most widely known fact is that we can constitute a relations between the invariants and features of a curve, as well as we could set the similar relations between one curve and another. This relation so interesting that up to now many papers have been written and still continues to be published. One of the commonly used exemplification revealing this relation is the Bertrand curve pair. We may cite some remarkable works drawing our attention. After analyzing the main properties of Bertrand curves, we see that geometers investigate these curve pairs in different spaces. While Babaarslan et al \cite{baba-2013} study main characters of these curves in Euclidean space, Ekmekci et al \cite{ekmekci-2001} make contributions as giving the necessary conditions to be a general helix in Lorentzian space. We may also mention that Yayli et al \cite{yayli-2017} define the Bertrand curve couples in a three dimensional Lie groups and give a relation between these curves and harmonic curvature function. It is another valuable study on curves that we can make a classification \cite{chen-1991}. After this paper we recognize that some curves can be called as biharmonic satisfying the condition that the Laplace image of mean curvature is equal to zero, while some of them provided the Laplace image of mean curvature is equal to a non-zero real constant times mean curvature, that is, 1-type of harmonic curve. Among many works we use as a tool only some of them: Senyurt et al \cite{cakir-2019}, they point out a method to classify a given curve by means of an another curve.
Kocayigit et al \cite{kocayigit-2011}, they study 1-type harmonic curves by using the mean curvature vector field of the curve itself and also we refer Senyurt et al \cite{cakir-2018}, they study biharmonic curves whose mean curvature vector field is the kernel of Laplacian. Throughout the present paper, we first take a regular curve and calculate the invariants of this curve. Then we name it as main(base) curve. By using the curvatures of the base curve we define the Bertrand partner of base curve. In order to determine the features of partner curve in terms of the curvatures of the base curve we apply the relations between a differentiable curve and its Bertrand mate. We give the conditions to be an harmonic curve (or being a 1-type harmonic) and also we derive the equations representing the Bertrand partner in terms of the main curve. Finally we give an example quite clearly support our allegations.\,\,\,
Frenet vectors are possible to give through the directional derivatives of these vectors, that is, traditional Frenet formulae are,\,\,see\cite{sabuncuoglu-2014}
\begin{eqnarray}\label{eq.(1.1)}
T'  =  \vartheta \kappa N,\,\,\,\,\,\,\,\,\,\,\,\,\,\,\,\,
N'  =  -\vartheta\kappa T + \vartheta\tau B,\,\,\,\,\,\,\,\,\,\,\,\,\,\,\,
B'  =  -\vartheta\tau N.
\end{eqnarray}
Given that $\alpha$ is a differentiable curve with the principal normal $N$ and $\beta$ is another differentiable curve. $\beta$ is called the Bertrand partner of $\alpha$ if they have the common principal normal at the corresponding points. Then $(\alpha, \beta)$ is called the Bertrand curve pair. It is obvious from this statement,\,\,see\cite{sabuncuoglu-2014}
\begin{eqnarray}\label{eq.(1.2)}
\beta(t) = \alpha(t) + \lambda(t) N(t),\,\,\,\,\,\,\lambda(t)\,\in\mathbb{R} .
\end{eqnarray}
\begin{figure}[h]
\centering
\includegraphics[height=0.15\textheight]{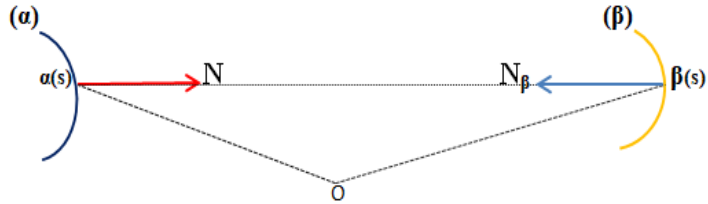}\centering\\
  \caption{($\alpha\,,\gamma$) Bertrand curve pairs}
\end{figure}
The relation between the Frenet vectors of $\alpha$ and  $\beta$ is
\begin{eqnarray}\label{eq.(1.3)}
&&T_{\beta}  =  cos\theta T + sin\theta B,\,\,\,\,\,
N_{\beta}  =  N,\,\,\,\,\,
B_{\beta}  =  -sin\theta T + cos\theta B,\\
&&\textrm{where}\,\,\,cos\theta = \,<T_{\beta}, T >\cdot\nonumber
\end{eqnarray}
Also the relation between the curvatures of $\alpha$ and  $ \beta$ is
\begin{eqnarray}\label{eq.(1.4)}
\kappa_{\beta}(s)  =  \frac{\lambda\kappa-sin^{2}\theta}{\lambda(1-\lambda\kappa)}\,\,,\,\,\,\,\,\,\,\,\,\,
\tau_{\beta}(s)  =  \frac{1}{\lambda^{2}\tau}sin^{2}\theta\,\cdot
\end{eqnarray}
The ordered pair $(\alpha, \beta)$ forms a Bertrand couple if and only if,\,\,see\,\cite{sabuncuoglu-2014}
\begin{eqnarray}\label{eq.(1.5)}
\lambda\kappa + \mu\tau = 1 \,\,\,\,\textrm{where}\,\,\,\lambda, \mu \in \mathbb{R}.
\end{eqnarray}
\,Given that\,\,$\alpha$\,is a unit speed curve with principal normal $N.$ We define the mean curvature vector of $\alpha$,\,\,\,see\,\cite{chen-1991}
\begin{eqnarray}\label{eq.(1.6)}
H = D_{\alpha'} \alpha' = \kappa N
\end{eqnarray}
where $D$ is the Levi-Civita connection.\,\,
According to this definition the mapping
\begin{eqnarray}\label{eq.(1.7)}
\Delta : \chi^{\bot} (\alpha(I)) \rightarrow \chi (\alpha(I))\,\,\,\,\,\textrm{such that}\,\,\,\,\,\Delta H = -D^{2}_{T} H
\end{eqnarray}
is called a Laplacian operator. When we phrase the normal bundle of
$\alpha$ \,by\, $\chi^{\bot}(\alpha(s))$\,, it follows that the normal connection  $D^{\perp}$ is
\begin{eqnarray}\label{eq.(1.8)}
D^{\bot}_{T} : \chi^{\bot}(\alpha(I)) \rightarrow \chi^{\bot}(\alpha(I))\,,\,\,\,\,\,\,\,\,
D^{\bot}_{T}X = D_{T}X - \big < D_{T}X , T \big > T
\end{eqnarray}
and the normal Laplacian operator $\Delta^{\bot}$ is
\begin{eqnarray}\label{eq.(1.9)}
\Delta^{\bot}_{T} X = -D^{\bot}_{T} D^{\bot}_{T}X,
 \,\,\,\forall X \in \chi^{\bot}(\alpha(I)).
\end{eqnarray}

\begin{theorem}\emph{\cite{chen-1991}}
\,Let $\alpha$ be a regular curve with the mean curvature \,$H.$ Then the following propositions hold,\,\,\\
\\
i)\,\,$\Delta H = 0$\,\,then\,\,$\alpha$\,\,is called biharmonic curve,\,\,\,ii) $\Delta H = \lambda H,$\,\,$\lambda \in \mathbb{R}$,\,\,provided\,\,$\alpha$\,\,is of 1-type harmonic.\\
\\
iii)\,\,$\Delta^{\perp} H = 0$\,\,then\,\,$\alpha$\,\,is called weak biharmonic curve\,,\,\,\,iv) $\Delta^{\perp} H = \lambda H,$\,\,$\lambda \in \mathbb{R}$,\,\,provided\,\,$\alpha$\,\,is of 1-type harmonic\,.
\end{theorem}
\begin{theorem}\emph{\cite{cakir-2018}}
Let a differentiable curve $\alpha $ be given. We have the following equation representing the curve $\alpha$ with respect to unit tangent $T$
$$
D^{3}_{T}T + \lambda_{2}D^{2}_{T}T + \lambda_{1}D_{T}T + \lambda_{0}T = 0
$$
with the coefficients $\lambda_{0}$, $\lambda_{1}$ and $\lambda_{2}$
\begin{eqnarray*}
&&\lambda_{0}  = \vartheta^{2} \kappa \tau (\frac{\kappa}{\tau})'\,,\,\,\,\,
\lambda_{1} =  \vartheta^{2} (\kappa^{2} + \tau^{2}) - \frac{\vartheta''}{\vartheta} - \frac{\kappa''}{\kappa} + (\frac{\vartheta'}{\vartheta} + \frac{\kappa'}{\kappa})(3 \frac{\vartheta'}{\vartheta}  +  \frac{\tau'}{\tau}) + 2(\frac{\kappa'}{\kappa})^{2}\,\,\,\\
&&\lambda_{2}  =  -(3\frac{\vartheta'}{\vartheta}+2\frac{\kappa'}{\kappa}+\frac{\tau'}{\tau})\,.
\end{eqnarray*}
\end{theorem}
\small\section{Discussions and Result}\label{S:intro}

\noindent Throughout the present paper we use the set $\{T,\,N,\,B,\,\kappa,\,\tau \}$ to express the Frenet apparatus of the main curve $\alpha$ and also the set $\{T_{\beta},\,N_{\beta},\,B_{\beta},\,\kappa_{\beta},\,\tau_{\beta} \}$ stands for the Frenet elements of the partner curve $\beta$ with the norm $\vartheta =\, \parallel\beta'(s)\parallel$.

\begin{theorem}
Let $(\alpha, \beta)$ be a Bertrand curve pair. Then the covariants derivatives of $\alpha$ with respect to $B$ is given
\begin{eqnarray}\label{denklem 4-4-1}
D_{B}T & = & \big(\frac{1-cos\theta}{sin\theta}\big)\kappa N, \nonumber \\
\nonumber \\
D_{B}N & = & -\big(\frac{1-cos\theta}{sin\theta}\big)\kappa T + \big(\frac{1-cos\theta}{sin\theta}\big)\tau B, \\
\nonumber \\
D_{B}B & = & -\big(\frac{1-cos\theta}{sin\theta}\big)\tau N.\nonumber
\end{eqnarray}
\begin{proof}
From eq.(\ref{eq.(1.1)}) and eq.(\ref{eq.(1.2)}) we have the norm \,\, $\vartheta =\, \parallel\beta'(s)\parallel = \tau\sqrt{\lambda^{2}+\mu^{2}}$ and it follows from the equalities (\ref{eq.(1.3)}) and (\ref{eq.(1.4)})
\begin{eqnarray}\label{denklem 4.2-7}
D_{T_{\beta}}T_{\beta} &=& \vartheta\kappa_{\beta}N_{\beta} \nonumber \\
\nonumber \\
& = & \vartheta\frac{\lambda\kappa - sin^{2}\theta}{\lambda(1-\lambda\kappa)}N
\end{eqnarray}
Employing the equalities (\ref{eq.(1.3)}) and (\ref{eq.(1.4)}) for the left side of this equality gives
\begin{eqnarray}\label{denklem 4.2-8}
D_{T_{\beta}}T_{\beta} &=& D_{(cos\theta T + sin\theta B)}(cos\theta T + sin\theta B) \nonumber \\
&=& cos\theta D_{T}(cos\theta T + sin\theta B) + sin\theta D_{B}(cos\theta T + sin\theta B)
\nonumber \\
&=& cos^{2}\theta\kappa N - cos\theta sin\theta \tau N + cos\theta sin\theta D_{B}T + sin^{2}\theta D_{B}B.
\end{eqnarray}
From (\ref{denklem 4.2-7}) and (\ref{denklem 4.2-8}) we get
\begin{eqnarray}\label{denklem 4.2-9}
  cos\theta sin\theta D_{B}T + sin^{2}\theta D_{B}B & = & \big(\vartheta\frac{\lambda\kappa - sin^{2}\theta}{\lambda(1-\lambda\kappa)}-\kappa cos^{2}\theta +\tau cos\theta sin\theta\big)N.
\end{eqnarray}
By the same way we write
\begin{eqnarray}\label{denklem 4.2-10.}
D_{T_{\beta}}N_{\beta} &=& -\vartheta\kappa_{\beta}T_{\beta}+\vartheta\tau_{\beta}B_{\beta} \nonumber \\
\nonumber \\
 & = & -\vartheta\frac{\lambda\kappa - sin^{2}\theta}{\lambda(1-\lambda\kappa)}\big(cos\theta T + sin \theta B\big) + \frac{\vartheta sin^{2} \theta}{\lambda^{2}\tau}\big(-sin \theta T + cos\theta B\big) \nonumber \\
\nonumber \\
 & = & \big(\frac{-\vartheta sin^{3}\theta}{\lambda^{2}\tau}-\vartheta\frac{\lambda\kappa - sin^{2}\theta}{\lambda(1-\lambda\kappa)}cos\theta \big)T + \big(\frac{\vartheta sin^{2}\theta cos\theta}{\lambda^{2}\tau}-\vartheta\frac{\lambda\kappa - sin^{2}\theta}{\lambda(1-\lambda\kappa)} sin \theta \big)B . \nonumber\\
\end{eqnarray}
From the equalities (\ref{eq.(1.3)}) and (\ref{eq.(1.4)}) left side of this equality is
\begin{eqnarray}\label{denklem 4.2-11.}
D_{T_{\beta}}N_{\beta} &=& D_{(cos\theta T+sin\theta B)}N \nonumber \\
\nonumber \\
&=& cos\theta D_{T}N + sin\theta D_{B}N \nonumber \\
\nonumber \\
&=& -\kappa cos\theta T + \tau cos\theta B + sin\theta D_{B}N.
\end{eqnarray}
From (\ref{denklem 4.2-10.}) and (\ref{denklem 4.2-11.}) we find out
\begin{eqnarray*}
D_{B}N & = & \Big(\kappa cot \theta-\vartheta\frac{\lambda\kappa-sin^{2}\theta}{\lambda(1-\lambda\kappa)}cot \theta - \frac{\vartheta sin^{2} \theta}{\lambda^{2}\tau}\Big)T \nonumber \\
\nonumber \\
&& +\Big(\frac{\vartheta sin \theta cos\theta}{\lambda^{2}\tau}-\tau cot \theta - \vartheta\frac{\lambda\kappa-sin^{2}\theta}{\lambda(1-\lambda\kappa)} \Big)B.
\end{eqnarray*}
\noindent From eq.(\ref{eq.(1.5)}), $\lambda \kappa + \mu \tau = 1$\,\,\,\,and we also have\,\,\,\,$\vartheta = \tau\sqrt{\lambda^{2}+\mu^{2}}$. Hence we conclude
\begin{eqnarray*}
D_{B}N & = & \big(\frac{cos\theta-1}{sin\theta}\kappa\big)T + \big(\frac{1-cos\theta}{sin\theta}\tau\big)B.
\end{eqnarray*}
Finally making use of the above methods
\begin{eqnarray}\label{denklem 4.2-13.}
D_{T_{\beta}}B_{\beta} &=& -\vartheta\tau_{\beta}N_{\beta} \nonumber \\
\nonumber \\
 & = & -\frac{\vartheta sin^{2}\theta}{\lambda^{2}\tau}N.
\end{eqnarray}
Left side of this equality is
\begin{eqnarray}\label{denklem 4.2-14.}
D_{T_{\beta}}B_{\beta} &=& D_{(cos \theta T + sin \theta B)}(-sin \theta T + cos\theta B)\nonumber \\
\nonumber \\
& = & cos\theta D_{T}(-sin \theta T + cos\theta B) + sin \theta D_{B}(-sin \theta T + cos\theta B) \nonumber \\
\nonumber \\
&=& -cos\theta sin \theta D_{T}T + cos^{2}\theta D_{T}B - sin^{2}\theta D_{B}T + sin \theta cos\theta D_{B}B \nonumber \\
\nonumber \\
&=& (-\kappa cos\theta sin \theta - \tau cos^{2}\theta )N - sin^{2}\theta D_{B}T + sin \theta cos\theta D_{B}B.
\end{eqnarray}
From (\ref{denklem 4.2-13.}) and (\ref{denklem 4.2-14.}) we get
\begin{eqnarray}\label{denklem 4.2-15}
-sin^{2}\theta D_{B}T + cos\theta sin\theta D_{B}B = (\kappa sin\theta cos\theta-\frac{\vartheta sin^{2}\theta}{\lambda^{2}\tau} + \tau cos^{2}\theta)N.
\end{eqnarray}
From the equalities\,\,(\ref{denklem 4.2-9}) and (\ref{denklem 4.2-15}) we can write the counterparts of\, $D_{B}T$ and $D_{B}B$ respectively,\\
\begin{eqnarray*}
D_{B}T = \Big( \frac{\vartheta sin^{2}\theta}{\lambda^{2}\tau} + \vartheta\frac{\lambda\kappa-sin^{2}\theta}{\lambda(1-\lambda\kappa)}cot\theta-\kappa cot\theta\Big)N
\end{eqnarray*}
and
\begin{eqnarray*}
D_{B}B = \Big(\vartheta\frac{\lambda\kappa-sin^{2}\theta}{\lambda(1-\lambda\kappa)} - \frac{\vartheta sin\theta cos\theta}{\lambda^{2}\tau} + \tau cot\theta\Big)N.\,
\end{eqnarray*}
\\
\,\,Again taking the eq.(\ref{eq.(1.5)}) into consideration\,\,\,and\,\,\,using \,$\vartheta = \tau\sqrt{\lambda^{2}+\mu^{2}}$\,\,\,we obtain
\begin{eqnarray*}
D_{B}T & = & \big(\frac{1-cos\theta}{sin\theta}\kappa\big)N
\end{eqnarray*}
and
\begin{eqnarray*}
D_{B}B & = & \big(\frac{cos\theta-1}{sin\theta}\tau\big)N.
\end{eqnarray*}

\end{proof}
\end{theorem}

\begin{theorem}
Let $(\alpha, \beta)$ be a Bertrand curve pair. Then we give the covariants derivatives of $\alpha$ with respect to normal Levi-Civita connection $D^{\perp}$ as
\begin{eqnarray}\label{eq.(2.10)}
D^{\perp}_{B}T  =  \big(\frac{1-cos\theta}{sin\theta}\big)\kappa N\,,\,\,\,\,\,\,\,\,
D^{\perp}_{B}N  =  -\big(\frac{1-cos\theta}{sin\theta}\big)\kappa T
\end{eqnarray}
\begin{proof}
It is obvious from eq.(\ref{denklem 4-4-1}) that
\begin{eqnarray*}
D_{B}T & = & \big(\frac{1-cos\theta}{sin\theta}\big)\kappa N,  \\
\\
D_{B}N & = & -\big(\frac{1-cos\theta}{sin\theta}\big)\kappa T + \big(\frac{1-cos\theta}{sin\theta}\big)\tau B, \\
 \\
D_{B}B & = & -\big(\frac{1-cos\theta}{sin\theta}\big)\tau N.
\end{eqnarray*}
Now it remains only to apply the eq.(\ref{eq.(1.8)}). When we make use of this relation we obtain
\begin{eqnarray*}
D^{\bot}_{B}X = D_{B}X - \big < D_{B}X , B \big > B\,\,\,\,\,\,\,\Longrightarrow\,\,\,\,\,\,\,D^{\perp}_{B}T  =  \big(\frac{1-cos\theta}{sin\theta}\big)\kappa N
\end{eqnarray*}
and
\begin{eqnarray*}
D^{\bot}_{B}X = D_{B}X - \big < D_{B}X , B \big > B\,\,\,\,\,\,\,\Longrightarrow\,\,\,\,\,\,\,D^{\perp}_{B}N  =  -\big(\frac{1-cos\theta}{sin\theta}\big)\kappa T.
\end{eqnarray*}
\end{proof}
\end{theorem}

\begin{theorem}
Let $(\alpha, \beta)$ be a Bertrand curve pair. Then the following statements are equivalent with respect to connection $D$.\\
1) $\beta$ partner curve is biharmonic if and only if
\begin{eqnarray}\label{denklem 4.2.21}
(3\lambda\kappa-sin^{2}\theta)\kappa' = 0\nonumber\\
\nonumber\\
(\lambda\kappa-sin^{2}\theta)(1-cos\theta)^{2}(\kappa^{2}+\tau^{2})-\lambda\kappa''sin^{2}\theta = 0\\
\nonumber\\
sin^{2}\theta\tau'-\lambda\kappa\tau'-2\lambda\kappa'\tau = 0 \nonumber
\end{eqnarray}

2) $\beta$ partner curve is 1-type of harmonic if and only if
\begin{eqnarray}\label{denklem 4.2.22}
(3\lambda\kappa-sin^{2}\theta)\kappa' = 0\nonumber\\
\nonumber\\
\frac{(\lambda\kappa-sin^{2}\theta)(1-cos\theta)^{2}(\kappa^{2}+\tau^{2})-\lambda\kappa''sin^{2}\theta}{ sin^{2}\theta} = \lambda\,(\lambda\kappa-sin^{2}\theta)\\
\nonumber\\
sin^{2}\theta\tau'-\lambda\kappa\tau'-2\lambda\kappa'\tau = 0, \nonumber\,\,\,\,\lambda \in \mathbb{R}.
\end{eqnarray}
\begin{proof}
From eq.(\ref{eq.(1.6)}) we write $\mathbb{H_{\beta}} = \vartheta\kappa_{\beta}N_{\beta}$\, and from the equalities (\ref{eq.(1.3)}) and (\ref{eq.(1.4)}) we find
\begin{eqnarray}\label{denklem 4.4.14}
\mathbb{H_{\beta}} = \vartheta\frac{\lambda\kappa - sin^{2}\theta}{\lambda(1-\lambda\kappa)}N= \big(\frac{\lambda\kappa-sin^{2}\theta}{\mu\,sin\theta}\big)\,N.
\end{eqnarray}
If we constrict the statements as
\begin{eqnarray*}
 \frac{\lambda\kappa-sin^{2}\theta}{\mu\,sin\theta} = a_{1}\,,\,\,\,\,\,\,\,\frac{1-cos\theta}{sin\theta}\kappa = \alpha_{1}\,,\,\,\,\,\,\,\,\frac{cos\theta-1}{sin\theta}\tau  = \alpha_{2}
\end{eqnarray*}
we can calculate the mean curvature more precisely. From eq.(\ref{eq.(1.1)}) and eq.(\ref{eq.(1.7)}) we evaluate the Laplace image of mean curvature $\mathbb{H_{\beta}}$ as
\begin{eqnarray*}
\Delta\mathbb{H_{\beta}} &=& -D^{2}_{B}(\mathbb{H_{\beta}}) = -D_{B}D_{B}(a_{1}N)\\
\\
&=& -D_{B}\big(a'_{1}N + a_{1}D_{B}N\big) \\
\\
&=& \big(a'_{1}\alpha_{1} + (a_{1}\alpha_{1})'\big)T + \big(a_{1}(\alpha^{2}_{1}+\alpha^{2}_{2}) - a''_{1} \big)N + \big(a'_{1}\alpha_{2} + (a_{1}\alpha_{2})' \big)B\\
\\
& = & \big(\frac{(3\lambda\kappa-sin^{2}\theta)\kappa'}{\mu(1+cos\theta)} \big)T + \big(\frac{(\lambda\kappa-sin^{2}\theta)(1-cos\theta)^{2}(\kappa^{2}+\tau^{2})-\lambda\kappa''sin^{2}\theta}{\mu sin^{3}\theta} \big)N \\
\\
&&+ \big(\frac{sin^{2}\theta\tau'-\lambda\kappa\tau'-2\lambda\kappa'\tau}{\mu(1+cos\theta)} \big)B
\end{eqnarray*}
If we consider the case $\Delta \mathbb{H}_{\beta} = 0$ we see that eq.(\ref{denklem 4.2.21}) is satisfied. Accordingly taking the condition
$\Delta \mathbb{H}_{\beta} = \lambda\mathbb{H}_{\beta}$ into consideration, another case eq.(\ref{denklem 4.2.22}) holds.
\end{proof}
\end{theorem}

\begin{theorem}
Let $(\alpha, \beta)$ be a Bertrand curve pair. Then the following statements are equivalent with respect to normal connection $D^{\perp}$.\\
1) $\beta$ partner curve is weak biharmonic if and only if
\begin{eqnarray}\label{denklem 4.2.24}
3\lambda\kappa\kappa'-sin^{2}\theta\,\kappa' = 0 \nonumber \\
\\
\frac{(\lambda\kappa-sin^{2}\theta)(1-cos\theta)^{2}\kappa^{2}-\lambda\kappa''sin^{2}\theta}{\mu sin^{3}\theta} = 0 \nonumber
\end{eqnarray}
2) $\beta$ partner curve is 1-type of harmonic if and only if
\begin{eqnarray}\label{denklem 4.2.25}
3\lambda\kappa\kappa'-sin^{2}\theta\,\kappa' = 0 \nonumber \\
\\
\frac{(\lambda\kappa-sin^{2}\theta)(1-cos\theta)^{2}\kappa^{2}-\lambda\kappa''sin^{2}\theta}{ sin^{2}\theta} =\lambda(\lambda\kappa-sin^{2}\theta) \nonumber
\end{eqnarray}
\begin{proof}
From eq.(\ref{eq.(1.6)}) we have $\mathbb{H_{\beta}} = \vartheta\kappa_{\beta}N_{\beta}$\, and from the equalities (\ref{eq.(1.3)}) and (\ref{eq.(1.4)}) we can easily figure out
\begin{eqnarray}\label{denklem 4.4.14}
\mathbb{H_{\beta}} = \vartheta\frac{\lambda\kappa - sin^{2}\theta}{\lambda(1-\lambda\kappa)}N = \big(\frac{\lambda\kappa-sin^{2}\theta}{\mu\,sin\theta}\big)N.
\end{eqnarray}
If we compress the phrases as
\begin{eqnarray*}
 \frac{\lambda\kappa-sin^{2}\theta}{\mu\,sin\theta} = a_{1}\,,\,\,\,\,\,\,\,\frac{1-cos\theta}{sin\theta}\kappa = \alpha_{1}\,,\,\,\,\,\,\,\,\frac{cos\theta-1}{sin\theta}\tau  = \alpha_{2}
\end{eqnarray*}
we may evaluate the mean curvature more precisely. From eq.(\ref{eq.(1.1)}) and eq.(\ref{eq.(1.9)}) we evaluate the Laplace image of mean curvature $\mathbb{H_{\beta}}$ as
\begin{eqnarray*}
\Delta^{\perp}\mathbb{H_{\beta}} & = & -D^{\perp}_{B}D^{\perp}_{B}(a_{1}N)\\
\\
& = & -D^{\perp}_{B}\big(a'_{1}N + a_{1}D^{\perp}_{B}N\big) \\
\\
& = & \big(\frac{3\lambda\kappa\kappa'-sin^{2}\theta\,\kappa'}{\mu(1+cos\theta)} \big)T + \big(\frac{(\lambda\kappa-sin^{2}\theta)(1-cos\theta)^{2}\kappa^{2}-\lambda\kappa''sin^{2}\theta}{\mu sin^{3}\theta} \big)N.
\end{eqnarray*}
If we consider the case $\Delta^{\perp} \mathbb{H}_{\beta} = 0$, then we see that the eq.(\ref{denklem 4.2.24}) is held.  In a similar case taking the phrase $\Delta^{\perp} \mathbb{H}_{\beta} = \lambda\mathbb{H}_{\beta}$ into consideration, we obtain that eq.(\ref{denklem 4.2.25}) is satisfied.

\end{proof}
\end{theorem}

\begin{corollary}
Let $(\alpha, \beta)$ be a Bertrand curve pair. Then the partner curve $\beta$ is of 1-type harmonic with respect to connection $D,$ provided that
$\alpha$ is a circular helix.\\
\begin{proof}  Given that
\,\,\,$\alpha(t) = \frac{1}{\sqrt{2}}(cost , sint , t)$\,\,\,be a circular helix. Then the Bertrand partner of this curve is
$$
\beta(t) = \frac{1}{\sqrt{2}}(cost , sint , t) + \lambda(-cost , -sint , 0)\,\,,\,\,\,\lambda \in \mathbb{R}\,.
$$
If we compute the Laplace image of $\mathbb{H_{\beta}}= \frac{\lambda\kappa-sin^{2}\theta}{\mu sin\theta}N$ from eq.
(\ref{denklem 4-4-1}) we figure out
$$
\Delta\mathbb{H_{\beta}} =\frac{\lambda\kappa-sin^{2}\theta}{\mu\,sin\theta}\,(\frac{cos\theta-1}{sin\theta})^{2}N.
$$
Considering the eq.(\ref{denklem 4.2.22}), we obtain $\beta$ is a 1-type of harmonic curve provided, \,\, $\lambda = (\frac{cos\theta-1}{sin\theta})^{2} \in \mathbb{R}$.

\end{proof}
\end{corollary}

\begin{corollary}
Let $(\alpha, \beta)$ be a Bertrand curve pair. Then the partner curve $\beta$ is of 1-type harmonic with respect to normal connection $D^{\perp},$ provided that $\alpha$ is a circular helix.\\
\begin{proof}  Given that
\,\,\,$\alpha(t) = \frac{1}{\sqrt{2}}(cost , sint , t)$\,\,\,be a circular helix. Then the Bertrand partner of this curve is
$$
\beta(t) = \frac{1}{\sqrt{2}}(cost , sint , t) + \lambda(-cost , -sint , 0)\,\,,\,\,\,\lambda \in \mathbb{R}\,.
$$
If we compute the normal Laplace image of\, $\mathbb{H_{\beta}}= \frac{\lambda\kappa-sin^{2}\theta}{\mu sin\theta}N$ from eq.(\ref{eq.(2.10)}) we figure out
$$
\Delta\mathbb{H_{\beta}} =\frac{\lambda\kappa-sin^{2}\theta}{\mu\,sin\theta}\,(\frac{cos\theta-1}{sin\theta})^{2}N.
$$
Considering the eq.(\ref{denklem 4.2.25}), we obtain $\beta$ is a 1-type of harmonic curve provided, \,\, $\lambda = \frac{1}{2}\big(\frac{cos\theta-1}{sin\theta}\big)^{2}\in \mathbb{R}$

\end{proof}
\end{corollary}

\begin{theorem}
Let $(\alpha, \beta)$ be a Bertrand curve pair. Then we give the equation representing the partner curve $\beta$ with respect to binormal $B$ as in
\begin{eqnarray}
D^{3}_{B}B - \Big(\frac{\kappa'}{\kappa}+2\frac{\tau'}{\tau} \Big) D^{2}_{B}B + \Big(\big(\frac{1-cos\theta}{sin\theta}\big)^{2}(\kappa^{2}+\tau^{2}) &+& 2(\frac{\tau'}{\tau})^{2}+\frac{\kappa'\tau'}{\kappa\tau}-\frac{\tau''}{\tau} \Big)D_{B}B \nonumber\\
\nonumber\\
&+& \Big(\kappa\tau\big(\frac{1-cos\theta}{sin\theta}\big)^{2}\big(\frac{\tau}{\kappa}\big)' \Big)B = 0.\nonumber\\
\end{eqnarray}

\begin{proof}
Let's put the eq.(\ref{denklem 4-4-1}) in another way as
\begin{eqnarray}\label{denklem 4.4.19}
D_{B}T = \alpha_{1}N\,,\,\,\,\,\,\,\,\,D_{B}N = -\alpha_{1}T - \alpha_{2}B\,,\,\,\,\,\,\,\,\,D_{B}B = \alpha_{2}N
\end{eqnarray}
with the coefficients \, $\alpha_{1},\,\,\alpha_{2}$
\begin{eqnarray*}
&&\alpha_{1} =\big(\frac{1-cos\theta}{sin\theta}\kappa\big)\,\,\,\,\,\,\,\,\,\,\,\textrm{and}\,\,\,\,\,\,\,\,\,\,\,
\alpha_{2} = \big(\frac{cos\theta-1}{sin\theta}\tau\big).
\end{eqnarray*}
It is open that we can derive
\begin{eqnarray}\label{denklem 4.4.20}
D_{B}B = \alpha_{2}N\,\,\,\Longrightarrow\,\,\,N = \frac{1}{\alpha_{2}}D_{B}B
\end{eqnarray}
and also
\begin{eqnarray}\label{denklem 4.4.21}
D_{B}T = \alpha_{1}N\,\,\,\,\textrm{and}\,\,\,\,D_{B}B = \alpha_{2}N\,\,\,\Longrightarrow\,\,\,D_{B}T = \frac{\alpha_{1}}{\alpha_{2}}D_{B}B.
\end{eqnarray}
Taking the consecutive derivatives of $D_{B}B = \alpha_{2}N$ with respect to $B$ we get
\begin{eqnarray*}
D_{B}B  =  \alpha_{2}N\,\,\,\,\Longrightarrow\,\,\,\,D_{B}(D_{B}B) & = & D_{B}(\alpha_{2}N)\\
\\
D^{2}_{B}B & = & -\alpha_{1}\alpha_{2}T + \alpha'_{2}N - \alpha^{2}_{2}B
\end{eqnarray*}
and
\begin{eqnarray*}
D_{B}(D^{2}_{B}B) & = & D_{B}(-\alpha_{1}\alpha_{2}T + \alpha'_{2}N - \alpha^{2}_{2}B)\\
\\
D^{3}_{B}B & = & (-\alpha_{1}\alpha_{2})'T - \alpha_{1}\alpha_{2}D_{B}T + \alpha''_{2}N + \alpha'_{2}D_{B}N - (\alpha^{2}_{2})'B - \alpha^{2}_{2}D_{B}B. \\
\end{eqnarray*}
Replacing the equivalents of $D_{B}N$, $N$ and $D_{B}T$ respectively from the equalities (\ref{denklem 4.4.19}), (\ref{denklem 4.4.20}) and (\ref{denklem 4.4.21}) we find
\begin{eqnarray}\label{denklem 4.4.22}
D^{3}_{B}B = \big(-(\alpha_{1}\alpha_{2})' - \alpha_{1}\alpha'_{2}\big)T  + \big(\frac{\alpha''_{2}}{\alpha_{2}}-\alpha^{2}_{1}-\alpha^{2}_{2}\big)D_{B}B - \big(\alpha_{2}\alpha'_{2} +( \alpha^{2}_{2})'\big)B.
\end{eqnarray}
Finally calculating the equivalent of \,$T$ we figure out
\begin{eqnarray*}
D_{B}N  =  -\alpha_{1}T - \alpha_{2}B \,\,\,\,\,\,\Longrightarrow\,\,\,\,\,\,T & = & -\frac{1}{\alpha_{1}}D_{B}N - \frac{\alpha_{2}}{\alpha_{1}}B \\
\\
 & = & \frac{-1}{\alpha_{1}\alpha_{2}}D^{2}_{B}B + \frac{\alpha'_{2}}{\alpha_{1}\alpha^{2}_{2}}D_{B}B - \frac{\alpha_{2}}{\alpha_{1}}B.
\end{eqnarray*}
Setting this phrase in eq.(\ref{denklem 4.4.22}) gives us
\begin{eqnarray*}
D^{3}_{B}B + \big((\alpha_{1}\alpha_{2})' + \alpha_{1}\alpha'_{2}\big)\big(\frac{-1}{\alpha_{1}\alpha_{2}}D^{2}_{B}B + \frac{\alpha'_{2}}{\alpha_{1}\alpha^{2}_{2}}D_{B}B - \frac{\alpha_{2}}{\alpha_{1}}B\big)&+&\big(\frac{-\alpha''_{2}}{\alpha_{2}} + \alpha^{2}_{1} + \alpha^{2}_{2}\big)D_{B}B \\
\\
&& + \big(\alpha_{2}\alpha'_{2} +( \alpha^{2}_{2})'\big)B. = 0
\end{eqnarray*}
If we rearrange the above equality we obtain
\begin{eqnarray*}
D^{3}_{B}B - \big(\frac{\alpha'_{1}}{\alpha_{1}} + 2\frac{\alpha'_{2}}{\alpha_{2}}\big) D^{2}_{B}B + \big(\alpha^{2}_{1} + \alpha^{2}_{2} + 2(\frac{\alpha'_{2}}{\alpha_{2}})^{2} + \frac{\alpha'_{1}\alpha'_{2}}{\alpha_{1}\alpha_{2}} - \frac{\alpha''_{2}}{\alpha_{2}}\big)D_{B}B + \big(\alpha'_{2}\alpha_{2} - \frac{\alpha'_{1}}{\alpha_{1}}\alpha^{2}_{2}\big)B = 0.
\end{eqnarray*}

\end{proof}
\end{theorem}

\begin{theorem}
Let $(\alpha, \beta)$ be a Bertrand curve pair. Then we give the following equations representing the partner curve $\beta$ with respect to normal connection $D^{\perp}$ as in
i) by means of unit tangent $T$
$$
\big(\frac{1-cos\theta}{sin\theta}\kappa\big)D^{\perp}_{B}D^{\perp}_{B}T - \big(\frac{1-cos\theta}{sin\theta}\kappa'\big)D^{\perp}_{B}T + \big(\frac{1-cos\theta}{sin\theta}\kappa\big)^{3}T = 0.
$$
ii) by means of principal normal $N$
$$
\big(\frac{1-cos\theta}{sin\theta}\kappa\big)D^{\perp}_{B}D^{\perp}_{B}N - \big(\frac{1-cos\theta}{sin\theta}\kappa'\big)D^{\perp}_{B}N + \big(\frac{1-cos\theta}{sin\theta}\kappa\big)^{3}N = 0.
$$
\begin{proof}
From eq.(\ref{eq.(2.10)}) we have $D^{\perp}_{B}T = \big(\frac{1-cos\theta}{sin\theta}\kappa\big)N$ and $D^{\perp}_{B}N = -\big(\frac{1-cos\theta}{sin\theta}\kappa\big)T$. Denoting the coefficients of $T\,$ and $\,N$ as\,\, $ \frac{1-cos\theta}{sin\theta}\kappa = \alpha_{1}$\,\,we can evaluate the derivatives more precisely
\begin{eqnarray*}
D^{\perp}_{B}T = \alpha_{1}N\,\,\,\Longrightarrow\,\,\,D^{\perp}_{B}\big(D^{\perp}_{B}T\big) & = & D^{\perp}_{B}\big(\alpha_{1}N\big) \\
& = & \alpha'_{1}N + \alpha_{1}D^{\perp}_{B}N \\
D^{\perp}_{B}D^{\perp}_{B}T& = & \frac{\alpha'_{1}}{\alpha_{1}}D^{\perp}_{B}T - \alpha^{2}_{1}T.
\end{eqnarray*}
If we rearrange this equality we obtain
$$
\alpha_{1}D^{\perp}_{B}D^{\perp}_{B}T - \alpha'_{1}D^{\perp}_{B}T + \alpha^{3}_{1}T = 0.
$$
In a similar way
\begin{eqnarray*}
D^{\perp}_{B}N = -\alpha_{1}T\,\,\,\Longrightarrow\,\,\,D^{\perp}_{B}\big(D^{\perp}_{B}N\big) & = & D^{\perp}_{B}\big(-\alpha_{1}T\big) \\
& = & -\alpha'_{1}T - \alpha_{1}D^{\perp}_{B}T \\
D^{\perp}_{B}D^{\perp}_{B}N& = & \frac{\alpha'_{1}}{\alpha_{1}}D^{\perp}_{B}N - \alpha^{2}_{1}N
\end{eqnarray*}
we obtain the differential equation by means of $N$ as
$$
\alpha_{1}D^{\perp}_{B}D^{\perp}_{B}N - \alpha'_{1}D^{\perp}_{B}N + \alpha^{3}_{1}N = 0.
$$

\end{proof}
\end{theorem}

\begin{example}Given that
\,\,\,$\alpha(t) = \frac{1}{\sqrt{2}}(cost , sint , t)$\,\,\,be a circular helix. Then the Bertrand partner of\, $\alpha$ is
$$
\beta(t) = \frac{1}{\sqrt{2}}(cost , sint , t) + \lambda(-cost , -sint , 0)\,\,,\,\,\,\lambda \in \mathbb{R}.
$$
After plain calculations we find out the differential equations characterizing the partner curve $\beta$ in terms of the main curve are \\
\\
i) According to eq.(\ref{denklem 4-4-1}) \\
$$
D^{3}_{B}B + (\frac{1-cos\theta}{sin\theta})^{2}D_{B}B = 0,
$$\\
ii) According to eq.(\ref{eq.(2.10)})
\begin{eqnarray*}
D^{\perp}_{B}D^{\perp}_{B}T + \big(\frac{cos\theta-1}{\sqrt{2}\,sin\theta}\big)^{2}T = 0
\end{eqnarray*}
and
\begin{eqnarray*}
D^{\perp}_{B}D^{\perp}_{B}N + \big(\frac{cos\theta-1}{\sqrt{2}\,sin\theta}\big)^{2}N = 0.
\end{eqnarray*}

\end{example}

\noindent \textbf{Conclusion:}
We give new Frenet formulae for the Bertrand partner curve by making use of relations between curvatures and a curve itself. Then referring these formulae we write the differential equations and sufficient conditions of harmonicity of the Bertrand partner curve through the main curve.  We hope this work inspire the geometers study in spaces other than that of Euclidean space.

\end{document}